\documentclass[11pt,a4paper,leqno]{amsart}
\input{avv.sty}
\title{Topics in Special Functions}
\author{G. D. Anderson$^1$, M. K. Vamanamurthy, and M. Vuorinen}
\date{}

\begin{document}

\begin{center}
{\tiny \tt FILE: om845.tex, 9 Oct. 2001}
\end{center}

\maketitle
\markboth{\textsc{G. D. ANDERSON, M. K. VAMANAMURTHY, AND
M. VUORINEN}}{\textsc{TOPICS IN SPECIAL FUNCTIONS}}

\begin{abstract}The authors survey recent results in special functions, particularly the gamma function and the Gaussian hypergeometric function.
\end{abstract}

%\vskip 0.5cm

%\tableofcontents

\footnotetext[1]{This paper is an outgrowth of an invited talk 
given by the first author at the 18th Rolf Nevanlinna Colloquium 
in Helsinki, Finland, in August 2000.}

\section{Introduction}\label{sect:intro}

Conformal invariants are powerful tools in the study of
quasiconformal mappings, and many of these have expressions in
terms of special functions.  For instance, the distortion results 
in geometric function theory, such as the quasiconformal
Schwarz Lemma, involve special functions.  A frequent task is
to simplify complicated inequalities, so as to clarify
the dependence on important parameters without
sacrificing sharpness. For these reasons we were led to study, as
an independent subject, various questions for special functions such
as monotoneity properties and majorants/minorants in terms of
rational functions. These new inequalities gave refined versions of 
some classical distortion theorems for quasiconformal maps.  
The classes of functions that occur include complete elliptic integrals, 
hypergeometric functions, and Euler's gamma function.
The main  part of our research is summarized in \cite{AVV5}.

In the later development most of our research has involved 
applications to geometric properties of quasiconformal maps.  
However, some of the questions concerning special functions, 
raised in \cite{AVV1}, \cite{AVV3}, and \cite{AVV5}, 
relate to special functions which are 
useful in geometric function theory in general, not just to 
quasiconformal maps.  In this survey our goal is to review 
the latest developments of the latter type, due to many 
authors \cite{A1}--\cite{A9}, 
\cite{AlQ1, AlQ2, AW, BPR1, BPR2, BPS, BP, EL, K1, K2, Ku}.

The methods used in these studies are based on classical analysis.
One of the technical tools is the {Monotone l'H\^opital's Rule},
stated in the next paragraph, which played an important
role in our work \cite{AVV4}--\cite{AVV5}.    
The authors discovered this result in \cite{AVV4}, unaware 
that it had been used earlier (without the name) as a technical tool
in differential geometry.  See \cite[p.\ 124, Lemma 3.1]{Ch} 
or \cite[p.\ 14]{AQVV} for relevant remarks.

\begin{lem}\label{lem:eka}
For $-\infty < a < b < \infty$,
let  $g$  and  $h$  be real-valued functions that are
continuous on $[a,b]$ and differentiable on $(a,b)$, with $h' \ne 0$ on
$(a,b)$.  If  $g'/h'$ is strictly increasing 
{\rm (}resp. decreasing{\rm )} on $(a,b)$, then the functions
\begin{equation*}
\frac{g(x)-g(a)}{h(x)-h(a)}\quad \text{and}\quad \frac{g(x)-g(b)}{h(x)-h(b)}
\end{equation*}
are also strictly increasing {\rm (}resp. decreasing{\rm )} on $(a,b)$.
\end{lem}

Graphing of the functions and computer experiments in
general played an important role in our work. For instance, the 
software that comes with the book \cite{AVV5} provides computer
programs for such experiments.

We begin this survey by discussing some recent results on the gamma 
function, including monotoneity and convexity properties and close 
approximations for the Euler-Mascheroni constant.  Hypergeometric 
functions have a very central role in this survey.  We give here 
a detailed proof of the so-called Elliott's identity for these 
functions following an outline suggested by Andrews, Askey, 
and Roy in \cite[p.\ 138]{AAR}. This identity contains, as a special case, 
the classical Legendre Relation and has been studied recently 
in \cite{KV} and \cite{BPSV}.  After this we discuss mean values, 
a topic related to complete elliptic integrals and their estimation,
and we present several sharp approximations for complete elliptic integrals.
We display inequalities 
for hypergeometric functions that generalize the Landen relation, 
and conclude the paper with a remark on recent work of
geometric mapping properties of hypergeometric functions as a function
of a complex argument.

This survey does not cover recent work on the applications of special 
functions to the change of distance under quasiconformal maps.  For this 
subject the interested reader may consult \cite{AVV5}.

\section{The $\Gamma$ and $\Psi$ functions}\label{sect:gammapsi}

Throughout this paper $\Gamma$ will denote Euler's gamma function, defined by
\begin{equation*}
\Gamma (z) = \int^{\infty}_0 e^{-t} t^{z-1}\,dt, \quad  \text{Re}\,z > 0,
\end{equation*}
and then continued analytically to the finite complex plane minus the set of
nonpositive integers. The recurrence formula $\Gamma(z+1)=z\,\Gamma(z)$
yields $\Gamma(n+1)=n!$ for any positive integer $n$.
We also use the fact that $\Gamma(1/2)=\sqrt{\pi}.$
The beta function is related to the gamma function by
$B(a,b)= \Gamma(a)\Gamma(b)/\Gamma(a+b)$.
The logarithmic derivative of the gamma function will be denoted,
as usual, by
\begin{equation*}
\Psi(z)\equiv\frac{d}{dz} \log \Gamma(z) = \frac{\Gamma'(z)}{\Gamma(z)}.
\end{equation*}

The Euler-Mascheroni constant $\gamma$ is defined as (see \cite{A2},
\cite{TY}, \cite{Y})
\begin{equation*}
\gamma \equiv \lim_{n\to\infty} D_n =0.57721 56649\dots ; \
D_n \equiv \sum_{k=1}^n\frac{1}{k}-\log n.
\end{equation*}
Then $\Psi(1)= \Gamma'(1)= - \gamma$ and $\Psi(1/2) =-\gamma -2 \log 2$.
For a survey of the gamma function see  \cite{G}, and for some 
inequalities for the gamma and psi functions see \cite{A1}.

\begin{ApprEM}
The convergence of the sequence
$D_n$ %\equiv \Sigma_{k=1}^n{1\over k} - \log n$
to $\gamma$ is very slow (the speed of convergence is studied
by Alzer \cite{A2}).  D. W. DeTemple \cite{De} studied 
a modified sequence which
converges faster and proved
\begin{equation*}
\frac{1}{24(n+1)^2} < R_n-\gamma < \frac{1}{24n^2},
\quad \text{where} \quad
R_n\equiv \sum_{k=1}^n\frac{1}{k} - \log\left(n+\frac{1}{2}\right).
\end{equation*}

%%%%%%%%%%%%%%%%%%%%%%%%%%%%%%%%%%%%%%%%%%%%%%%%%%%
Now let
$$
h(n)=R_n-\gamma,\ \ H(n)=n^2h(n),\ n\geqslant 1.
$$
Since $\Psi(n) =-\gamma-1/n+ \sum_{k=1}^n{1/k},$ we see that
$$H(n)=(R_n-\gamma)n^2 = (\Psi(n)+1/n -\log(n+ \frac{1}{2}))n^2.$$
Some computer experiments led M. Vuorinen to conjecture that
$H(n)$ increases on the interval $[1,\infty)$
from $H(1) = -\gamma + 1 - \log (3/2) = 0.0173\ldots$ to
$1/24 = 0.0416\ldots$.  E. A. Karatsuba proved in \cite{K1}
that for all integers $n \geqslant 1, H(n) < H(n+1),$
 by clever use of Stirling's formula
and Fourier series.  Moreover, using the relation 
$\gamma = 1-\Gamma'(2)$ she obtained, for $ k \ge 1$,
% OLD FORMULA: $$\gamma=1-(\log k)
% \sum_{r=1}^{12k+1} \frac{(-1)^{r-1}k^{r+1}}{(r-1)!(r+1)}+
% \sum_{r=1}^{12k+1}\frac{(-1)^{r-1}k^{r+1}}{(r-1)!(r+1)^2}+O(2^{-k}),
% $$
$$    % New formula (not in the preprint;only in paper!!!)
%- {2 \over (12k)!} -2 k^2 e^{-k}
-c_k \le 
\gamma -1 + (\log k)
 \sum_{r=1}^{12k+1} d(k,r) 
- \sum_{r=1}^{12k+1}\frac{d(k,r)}{r+1} 
\le c_k ,\\
$$
where
$$
c_k=\frac{2}{(12k)!} +2 k^2 e^{-k},\quad 
d(k,r)= (-1)^{r-1} \frac{k^{r+1}}{(r-1)! (r+1)},
$$
giving exponential convergence.
Some computer experiments also seem to indicate that 
$((n+1)/n)^2 H(n)$ is a decreasing convex function.
\end{ApprEM}

\begin{GammaFG}
Formulas for geometric objects, such as volumes of solids and arc lengths of
curves, often involve special functions.  For example, if $\Omega_n$
denotes the volume of the
unit ball $B^n=\{x:|x|< 1\}$ in $\R^n$, and if $\omega_{n-1}$ denotes
the $(n-1)$-dimensional surface area of the unit sphere $
S^{n-1} = \{x:|x|=1\}$, $n\geqslant 2$, then
\begin{equation*}
\Omega_n = \frac{\pi^{n/2}}{\Gamma ((n/2)+1)} \, ; \
{\omega_{n-1}} = n \Omega_n .
\end{equation*}
It is well known that for $n \geqslant 7$ both $\Omega_n$ and $\omega_n$
decrease to $0$ (cf. \cite[2.28]{AVV5}).
However, neither $\Omega_n$ nor $\omega_n$ is monotone for
$n$ on $[2,\infty)$. On the other hand,
$\Omega_n^{1/(n \log n)}$ decreases to $e^{-1/2}$ as $n \to \infty$ 
\cite[Lemma 2.40(2)]{AVV1}.

Recently H. Alzer \cite{A4} has obtained the best possible constants 
$a,\, b,\, A,\, B,$ \linebreak $\alpha,\, \beta$ such that
\begin{align*}
a\,\Omega^{\frac{n}{n+1}}_{n+1} &\leqslant \Omega_n \leqslant 
b\,\Omega^{\frac{n}{n+1}}_{n+1},\\
\sqrt{\frac{n+A}{2 \pi}} &\leqslant \frac{\Omega_{n-1}}{\Omega_{n}} 
\leqslant \sqrt{\frac{n+B}{2 \pi}},\\
\left(1 + \frac{1}{n} \right)^{\alpha} &\leqslant 
\frac{\Omega^2_n}{\Omega_{n-1} \Omega_{n+1}} \leqslant 
\left( 1 + \frac{1}{n} \right)^{\beta}
\end{align*}
for all integers $n \geqslant 1$. He showed that 
$a = 2/\sqrt{\pi} = 1.12837\ldots$, $b = \sqrt{e} = 1.64872\dots$, 
$A = 1/2$, $B= \pi/2 - 1 = 0.57079\dots$, 
$\alpha = 2 - (\log \pi)/\log 2 = 0.34850\dots$, $\beta = 1/2$.
For some related results, see \cite{KlR}.
\end{GammaFG}

\begin{MonProp}
In \cite{AnQ} it is proved that the function
\begin{equation}\label{eq:eka}
f(x)\equiv \frac{\log \Gamma (x+1)}{ x\, \log x}
\end{equation}
is strictly
increasing from $(1,\infty )$ onto $(1-\gamma,1)$.
In particular, for $x\in (1,\infty )$,
\begin{equation}\label{eq:toka}
x^{(1-\gamma )x-1} < \Gamma (x) < x^{x-1}.
\end{equation}
The proof required the following two technical
lemmas, among others:
\end{MonProp}

\begin{lem}\label{lem:toka}
The function
\begin{equation*}
g(x)\equiv \sum_{n=1}^{\infty }
\frac{n-x}{(n+x)^3}
\end{equation*}
is positive for $x\in [1,4)$.
\end{lem}

\begin{lem}\label{lem:kolmas}
The function
\begin{equation}\label{eq:hx}
h(x)\equiv x^2\,\Psi'(1+x) -x\,\Psi (1+x) + \log \Gamma (1+x)
\end{equation}
is positive for all
$x\in [1,\infty )$.
\end{lem}

It was conjectured in \cite{AnQ} that the function $f$ in 
\eqref{eq:eka} is concave on $(1,\infty )$.

\begin{HA}
Horst Alzer \cite{A2} has given an elegant proof of the
monotoneity of the function $f$ in \eqref{eq:eka}
by using the Monotone l'H\^opital's Rule and the convolution theorem for 
Laplace transforms. In a later paper \cite{A3} he has improved 
the estimates in \eqref{eq:toka} to
\begin{equation}\label{eq:2-10}
x^{\alpha(x-1) - \gamma} < \Gamma(x) < x^{\beta(x-1) - \gamma},
\quad x \in (0,1),
\end{equation}
where $\alpha \equiv 1 - \gamma = 0.42278\dots$,
$\beta \equiv \frac{1}{2}\left(\pi^2/6 - \gamma\right) =
0.53385\dots$ are best possible.
If $x \in (1, \infty)$, he also showed that \eqref{eq:2-10} holds
with best constants
$\alpha \equiv \frac{1}{2} \left(\pi^2/6 - \gamma\right)=0.53385\dots$, 
$\beta \equiv 1$. 
\end{HA}

\begin{EL}
Elbert and Laforgia  \cite{EL} have shown that the function 
$g$ in Lemma \ref{lem:toka}
is positive for all $x>-1$.  They used this result to prove that
the function $h$ in Lemma \ref{lem:kolmas} is strictly decreasing 
from $(-1,0]$ onto $[0,\infty)$ and strictly increasing from 
$[0,\infty)$  onto $[0,\infty)$. They also showed that 
$f'' < 0$ for $x>1$, thus proving the Anderson-Qiu
conjecture \cite{AnQ}, where $f$ is as in \eqref{eq:eka}.
\end{EL}

\begin{BP}
Berg and Pedersen \cite{BP} have shown that the function $f$ in 
\eqref{eq:eka} is not only strictly increasing from $(0,\infty)$ 
onto $(0,1)$, but is even a (nonconstant) so-called \emph{Bernstein function}.
That is, $f > 0$ and $f'$ is completely monotonic, i.e., 
$f' > 0$, $f'' < 0$, $f'''> 0$, \dots. 
In particular, the function $f$ is strictly
increasing and strictly concave on $(0,\infty )$.

In fact, they have proved the stronger result that $1/f$ is a 
Stieltjes transform, that is, can be written in the form
$$ \frac{1}{f(x)} = c + \int^{\infty}_0 \frac{d \sigma(t)}{x+t}, \quad x > 0,$$
where the constant $c$ is non-negative and $\sigma$ is a non-negative 
measure on $[0,\infty)$ satisfying
$$ \int^{\infty}_0 \frac{d \sigma(t)}{1 + t} < \infty.$$
In particular, for $1/f$ they have shown by using Stirling's formula 
that $c = 1$.  Also they have obtained $d \sigma(t) = H(t) dt$, 
where $H$ is the continuous density
$$ H(t) = \left\{ \begin{array}{ll}
t \dfrac{\log | \Gamma(1-t) | + (k-1) \log t}{(\log |\Gamma (1-t)|)^2 + 
(k-1)^2 \pi^2}, &  t \in (k-1,k), k=1,2,\dots, \\[.4cm]
0 \hspace{1.9in}, & t = 1,2,\dots .  \end{array}\right. $$
Here log denotes the usual natural logarithm.  The density $H(t)$ tends 
to $1/\gamma$ as $t$ tends to $0$, and $\sigma$ has no mass at $0$.
\end{BP}

\begin{Rama}
In ``The Lost Notebook and Other Unpublished Papers'' of Ramanujan [Ra1], 
the Indian mathematical genius, appears the following record:
\begin{equation*}
``\Gamma (1+x) =\sqrt{\pi} \Big(\frac{x}{e}\Big)^x
\Big\{8x^3+4x^2+x+\frac{\theta_x}{30}\Big\}^{1/6},
\end{equation*}
where $\theta_x$ is a positive proper fraction
\begin{align*}
&\theta_0 =\frac{30}{\pi^3}=.9675\\
&\theta_{1/12} =.8071\quad \theta_{7/12}=.3058\\
&\theta_{2/12} =.6160\quad \theta_{8/12}=.3014\\
&\theta_{3/12} =.4867\quad \theta_{9/12}=.3041\\
&\theta_{4/12} =.4029\quad \theta_{10/12}=.3118\\
&\theta_{5/12} =.3509\quad \theta_{11/12}=.3227\\
&\theta_{6/12} =.3207\quad \theta_1=.3359\\
&\theta_{\infty }=1.\text{''}
\end{align*}
\end{Rama}

Of course, the values in the above table, except $\theta_{\infty}$, 
are irrational and hence the decimals should be nonterminating as 
well as nonrecurring.  The record stated above has been 
the subject of intense investigations and is reviewed 
in \cite{BCK}, page 48 (Question 754).  This note of Ramanujan 
led the authors of \cite{AVV5} to make the following conjecture.

\begin{conj}\label{conj:eka}
Let
\begin{equation*}
G(x) = (e/x)^x\Gamma (1+x)/\sqrt{\pi}
\end{equation*}
and
\begin{equation*}
H(x)=G(x)^6-8x^3-4x^2-x = \frac{\theta_x}{30}.
\end{equation*}
Then $H$ is increasing from $(1,\infty)$ into $(1/100,1/30)$ 
\cite[p.\ 476]{AVV5}.
\end{conj}
\begin{Kara}
In a nice piece of work, E. A. Karatsuba \cite{K2} has proved the 
above conjecture. She did this by representing the function $H(x)$ 
as an integral for which she was able to  find an asymptotic 
development. Her work also led to an interesting asymptotic formula 
for the gamma function:
\begin{equation*}
\begin{split}
&\Gamma (x+1)=\sqrt{\pi }\Big(\frac{x}{e}\Big)^x \Big(8x^3+4x^2+x+
\frac{1}{30}-\frac{11}{240x} + \frac{79}{3360x^2}  + 
\frac{3539}{201600x^3}\\
&-\frac{9511}{403200x^4}-\frac{10051}{716800x^5}+
\frac{47474887}{1277337600x^6}
+\frac{a_7}{x^7}+\cdots+\frac{a_n}{x^n} +\Delta_{n+1}(x)\Big)^{1/6},
\end{split}
\end{equation*}
\medskip
where $\Delta_{n+1}(x)=O(\frac{1}{x^{n+1}})$, as $x \to \infty$, 
and where each $a_k$ is given explicitly in terms of the Bernoulli numbers.
\end{Kara}

\medskip

\section{Hypergeometric functions}\label{sect:hypfun}

Given complex numbers
$a$, $b$, and $c$ with $c\neq 0,-1,-2, \dots $,
the \emph{Gaussian hypergeometric function} is the analytic
continuation to the slit plane $\C\setminus [1,\infty)$ of
\begin{equation}\label{eq:kolmas}
F(a,b;c;z)\!= \!{}_2 F_1(a,b;c;z)\! \equiv\!
\sum_{n=0}^{\infty} \frac{(a,n)(b,n)}{(c,n)} \frac{z^n}{n!}, 
\quad |z|<1.
\end{equation}
Here $(a,0)=1$ for $a\neq 0$, and $(a,n)$
is the \emph{shifted factorial function}
\begin{equation*}
(a,n)\equiv a(a+1)(a+2) \cdots (a+n-1)
\end{equation*}
for $n=1,2,3,\ldots $.

The hypergeometric function $w = F(a,b;c;z)$ in \eqref{eq:kolmas} 
has the simple differentiation formula
\begin{equation}\label{eq:neljas}
\frac{d}{dz}\,F(a,b;c;z)=\frac{ab}{c}\,F(a+1,b+1;c+1;z).
\end{equation}

The behavior of the hypergeometric function near $z=1$ in the 
three cases $ a+b < c$,
$a+b = c$, and $ a+b > c, ~ a,b,c > 0$, is given by
\begin{equation}\label{eq:viides}
\begin{cases}
F(a,b;c;1) = \frac{\Gamma(c) \Gamma(c-a-b)}{\Gamma(c-a) \Gamma(c-b)}, 
\ a+b < c,\\
B(a,b)F(a,b;a+b;z)+\log(1-z) \\
\qquad \qquad \quad = R(a,b)+ O((1-z)\log(1-z)),\\
F(a,b;c;z) = (1-z)^{c -a -b} F(c-a,c-b;c;z),\ c < a+b,
\end{cases}
\end{equation}
where $R(a,b) = -2 \gamma - \Psi(a) - \Psi(b)$,
$R(a) \equiv R(a,1-a)$, $R(\frac{1}{2}) = \log 16,$ and where log 
denotes the principal branch of the complex logarithm.  
The above asymptotic formula for the \emph{zero-balanced} case $a+b=c$ 
is due to Ramanujan (see \cite{As}, \cite{Be1}). This formula is  
implied by \cite[15.3.10]{AS}.

The asymptotic formula \eqref{eq:viides} gives a precise description 
of the behavior of
the function $F(a,b;a+b;z)$ near the logarithmic singularity $z=1$.
This singularity can be removed by an exponential change of variables 
and the transformed function will be nearly linear.

\begin{thm}\label{thm:eka}
\cite{AQVV}
For  $a,b > 0$, let $k(x) = F(a,b;a+b;1-e^{-x})$, $x > 0$.
Then $k$ is an increasing and convex function with $k'((0, \infty)) =
(ab/(a+b)$, $\Gamma(a+b)/(\Gamma(a) \Gamma(b)))$.
\end{thm}

\begin{thm}\label{thm:toka}
\cite{AQVV}
Given $a,b > 0$, and  $a+b > c$, $d \equiv a+b-c$, the function
$\ell (x) = F(a,b;c;1-(1+x)^{-1/d})$, $ x > 0$, is increasing and
convex, with $\ell '((0,\infty)) = (ab/(cd)$,
$\Gamma(c)\Gamma(d)/(\Gamma(a)\Gamma(b)))$.
\end{thm}

\begin{GaussCont}
The six functions $F(a\pm 1, b;c;z)$, $F(a, b\pm 1; c;z)$, 
$F(a, b;c \pm 1;z)$ are
called \emph{contiguous} to $F(a,b;c;z)$. Gauss discovered 15 relations
between $F(a,b;c;z)$ and pairs of its contiguous functions 
\cite[15.2.10--15.2.27]{AS}, \cite[Section 33]{R2}. 
If we apply these relations to the differentiation formula \eqref{eq:neljas}, we obtain 
the following useful formulas.
\end{GaussCont}

\begin{thm}\label{thm:kolmas}
For $a,b, c > 0$, $z \in (0,1)$,  let $u = u(z) =
F(a-1,b;c;z)$, $v = v(z) = F(a,b;c;z)$,
$u_1 = u(1-z)$, $v_1 = v(1-z)$.  Then
\begin{align}
z\frac{du}{dz} &= (a -1) (v - u),\label{eq:kuudes}\\
z(1-z)\frac{dv}{dz} &= (c -a) u + (a-c+bz)  v,\label{eq:seits}
\end{align}
and
\begin{equation}\label{eq:kahd}
\frac{ab}{c} z(1-z) F(a+1,b+1;c+1;z) =(c-a) u
+(a-c+bz)v.
\end{equation}
Furthermore,
\begin{equation}\label{eq:kahd2}
%\begin{cases}
 z(1\!-\!z) \dfrac{d}{dz}\bigl(uv_1 \!+\! u_1 v \!-\! vv_1 \bigr)\! =\!  (1\! -\! a\! -\! b) \bigl[(1\!-\!z)u v_1 \!-\! z u_1 v\! -\! (1\!-\!2z)v v_1) \bigr].
%\end{cases}
\end{equation}
\end{thm}

Formulas \eqref{eq:kuudes}-\eqref{eq:kahd} in Theorem \ref{thm:kolmas} are
well known. See, for example, \cite[2.5.8]{AAR}.  On the other hand, formula
\eqref{eq:kahd2}, which follows from \eqref{eq:kuudes}-\eqref{eq:seits}
is first proved in  \cite[3.13 (4)]{AQVV}.

Note that the formula
\begin{equation}\label{eq:yhd}
z(1-z) \frac{dF}{dz} = (c-b)F(a,b-1;c;z)+(b-c+az)F(a,b;c;z)
\end{equation}
follows from  \eqref{eq:seits} if we use 
the symmetry property $F(a,b;c;z) = F(b,a;c;z)$. 

\bigskip

\begin{cor}\label{cor:3.13}
With the notation of Theorem {\rm \ref{thm:kolmas}}, 
if $a \in (0,1),~ b= 1-a < c,$ then
$$ uv_1 + u_1 v - vv_1 = u(1) = 
\frac{(\Gamma (c))^2} {\Gamma (c+ a-1) \Gamma(c-a+1)}.$$
\end{cor}

\medskip

\medskip

\section{Hypergeometric differential equation
}\label{sect:HDE}
\medskip

The function $F(a,b;c;z)$ satisfies
the hypergeometric differential equation
\begin{equation}\label{eq:kymm}
z(1-z)w''+[c-(a+b+1)z]w'-abw=0.
\end{equation}
Kummer discovered solutions of \eqref{eq:kymm} in various domains, 
obtaining 24 in all; for a complete list of his solutions 
see \cite[pp.\ 174, 175]{R2}.

\begin{lem}\label{lem:hde1}
{\rm (1)} If $2c=a+b+1$ then both  $F(a,b;c;z)$ and  $F(a,b;c;1-z)$
satisfy {\rm (\ref{eq:kymm})} in the lens-shaped region 
$\{z:0<|z|<1,\ 0<|1-z|<1\}.$\linebreak
{\rm (2)} If $2c=a+b+1$ then both $F(a,b;c;z^2)$ and $F(a,b;c;1-z^2)$
satisfy the differential equation
\begin{equation}\label{eq:kymm2} 
z(1-z^2)w''+[2c-1-(2a+2b+1)z^2]w'-4abz w=0
\end{equation} 
in the common part of the disk $\{z:|z|<1\}$ and the lemniscate 
$\{z:|1-z^2|<1\}$.
\end{lem}

       {\bf Proof.} 
By Kummer (cf. \cite[pp. 174-177]{R2}), the functions 
$F(a,b;c;z)$ and $F(a,b;a+b+1-c;1-z)$ are solutions of (\ref{eq:kymm}) 
in $\{z:0<|z|<1\}$ and $\{z:0<|1-z|<1\}$, respectively.  
But $a+b+1-c = c$ under the stated hypotheses.  The result (2) 
follows from result (1) by the chain rule.\ \ $\square$

\medskip

\begin{lem}\label{lem:hde2}
The function $F(a,b;c;\sqrt{1-z^2})$ satisfies the differential equation
$$
Z^3(1-Z)zw''-\{Z(1-Z)+[c-(a+b+1)Z]Zz^2\}w'-abz^3w=0,
$$
in the subregion of the right half-plane bounded by the lemniscate  $r^2 = 2\cos (2\vartheta )$, $-\pi /4 \le \vartheta \le \pi /4$, $z=re^{i\vartheta}$.  Here $Z=\sqrt{1-z^2}$, where the square root indicates the principal branch.
\end{lem}

{\bf Proof.}
From (4.1), the differential equation for $w=F(a,b;c;t)$ is given by 
$$t(1-t)\frac{d^2w}{dt^2}+[c-(a+b+1)t]\frac{dw}{dt}-abw=0.$$
Now put $t=\sqrt{1-z^2}.$  Then
$$
\frac{dz}{dt}=-\frac{t}{z},\ \frac{dt}{dz}=-\frac{z}{t},\ \frac{d^2t}{dz^2}=-\frac{1}{t^3}
$$
and
$$
\frac{dw}{dt}=-\frac{t}{z}\frac{dw}{dz},\ \frac{d^2w}{dt^2}=\frac{t^2}{z^2}\frac{d^2w}{dz^2}-\frac{1}{z^3}\frac{dw}{dz}.
$$
So
$$
t(1-t)\Big[\frac{t^2}{z^2}w''-\frac{1}{z^3}w'\Big]+\Big[c-(a+b+1)t\Big]
\Big(-\frac{t}{z}\Big)w'-abw=0.
$$
Multiplying through by $z^3$ and replacing $t$ by $Z\equiv \sqrt{1-z^2}$ gives the result.$ \ \square$

\bigskip

If $w_1$ and $w_2$ are two solutions of a second order differential 
equation, then their \emph{Wronskian} is defined to be 
$W(w_1,w_2)\equiv w_1w_2'-w_2w_1'$.

\begin{lem}\label{lem:linind}
\cite[Lemma 3.2.6]{AAR}  
If $w_1$ and $w_2$ are two linearly independent solutions of 
\eqref{eq:kymm}, then
\begin{equation*}
W(z)=W(w_1,w_2)(z)=\frac{A}{z^c(1-z)^{a+b-c+1}},
\end{equation*}
where $A$ is a constant.
\end{lem}

(Note the misprint in \cite[(3.10)]{AAR}, where the coefficient 
$x(1-x)$ is missing from the first term.)
 
\begin{lem}\label{lem:vv1}
If $2c=a+b+1$ then, in the notation of Theorem {\rm \ref{thm:kolmas}}, 
\begin{equation}\label{eq:duren}
 (c-a)(uv_1 + u_1v) + (a-1)vv_1 = A \cdot z^{1-c}(1-z)^{1-c}.
\end{equation}
\end{lem}

{\bf Proof.} If $2c=a+b+1$ then by Lemma \ref{lem:hde1}(1), both $v(z)$
and $v(1-z)$ are solutions of (\ref{eq:kymm}).  Since
$W(z) = W(v_1,v)(z)= v'(z)v_1(z) -v(z)v_1'(z)$, we have 
\begin{eqnarray*}
 z(1-z)W(z) & =& z(1-z)(v'v_1 -vv_1')\\
 &=& (c-a)(uv_1 + u_1v) + (2a+b-2c)vv_1 \\
 &=& (c-a)(uv_1 + u_1v) + (a-1) vv_1.
\end{eqnarray*}

Next, since $2c=a+b+1$, Lemma \ref{lem:linind} shows that 
$z^c(1-z)^c W(z) = A$, and the result follows.\ \ $\square$
\medskip

Note that in the particular case $c=1, a=b= \frac{1}{2}$ 
the right side of (\ref{eq:duren}) is constant and the result is similar to
Corollary \ref{cor:3.13}. This particular case is
Legendre's Relation (\ref{eq:legendre}), and this proof of it is due to
Duren \cite{Du}.

\medskip
%%%%%%%%%%%%%%%%%
%%%%%%%%%%%%%%%%%
%%%%%%%%%%%%%%%%%

\begin{lem}\label{lem:vv11}
 If $a, b > 0, c \ge 1,$ and $2c = a+b+1,$ 
then the constant $A$ in Lemma {\rm \ref{lem:vv1}} is given by
$A = (\Gamma(c))^2 / (\Gamma(a) \Gamma(b)).$ 
In particular, if $c = 1$ then Lemma {\rm \ref{lem:vv1}} reduces 
to Legendre's Relation {\rm (\ref{eq:legkaea})} for generalized elliptic integrals.
\end{lem}

{\bf Proof.} The idea of the proof is to replace the possibly
unbounded hypergeometric functions in formula (\ref{eq:duren}) by bounded
or simpler ones. Therefore we consider three cases.

\medskip

{\emph Case} (1): $c \ge 2.$ Now $ a+b \ge c+1 \ge 3.$
 By (\ref{eq:viides}) or  \cite[15.3.3]{AS}, we have
$$
u(z) = (1-z)^{2-c} F(c+1-a, c-b; c; z), \quad
u_1(z) = z^{2-c} F(c+1-a, c-b; c; 1-z),
$$
$$
v(z) = (1-z)^{1-c} F(c-a, c-b; c; z),\quad
v_1(z) = z^{1-c} F(c-a, c-b; c; 1-z).
$$
Hence
$$
A= (c-a) [(1-z)F(c+1-a,c-b;c; z) F(c-a, c-b; c; 1-z)$$
$$ 
+ z F(c+1-a,c-b;c;1-z)F(c-a,c-b;c;z)]
$$
$$
    \quad  + (a-1) F(c-a,c-b;c;z)F(c-a,c-b;c;1-z). 
$$
Now, since $a+b - c = c-1,$ letting $z \to 0,$ from 
(\ref{eq:viides}) we get
\begin{eqnarray*}
A& = &(c-a) \frac{\Gamma(c) \Gamma(c-1)}{\Gamma(a)\Gamma(b))} + 
(a-1) \frac{\Gamma(c) \Gamma(c-1)}{\Gamma(a)\Gamma(b)}\\
&=& (c-1) \frac{\Gamma(c) \Gamma(c-1)}{\Gamma(a)\Gamma(b)} = 
\frac{(\Gamma(c))^2}{\Gamma(a)\Gamma(b)},
\end{eqnarray*}
as claimed.

\medskip

{\emph Case} (2): $1 < c < 2$. Now, $1 < c < a+b < c+1 < 3.$ Then
$$
A = (c-a) [ (1-z)^{c-1} u(z)F(c-a,c-b;c;1-z) + z^{c-1} u_1(z)F(c-a,c-b;c;z)]
  $$
$$
        + (a-1)F(c-a,c-b;c;z)F(c-a,c-b;c;1-z).
$$
 Now letting $z \to 0,$ from (\ref{eq:viides}), 
as in Case (1), we get
\begin{eqnarray*}
A & = &(c-a) \frac{\Gamma(c) \Gamma(c-1)}{\Gamma(a)\Gamma(b)} + 
(a-1) \frac{\Gamma(c) \Gamma(c-1)}{\Gamma(a)\Gamma(b)} \\
& = & (c-1) \frac{\Gamma(c) \Gamma(c-1)}{\Gamma(a)\Gamma(b)} =
 \frac{(\Gamma(c))^2}{\Gamma(a)\Gamma(b)},
\end{eqnarray*}
 as claimed.

\medskip

{\emph Case} (3): $c = 1.$ Now $a+b = 1.$ Then
$$A = (1-a) [ u(z)v_1(z) + u_1(z)v(z) - v(z)v_1(z) ]
$$
$$                
\qquad  = (1-a) u_1(z)v(z) +(1-a)v_1(z)[u(z)-v(z)].
$$
From \cite[ Ex. 21(4), p.71]{R1} we have
\begin{eqnarray*}
u(z) - v(z) & = & F(a-1,b;c;z) - F(a,b;c;z)\\
& = & \frac{c-b}{c}zF(a,b;c+1;z) - zF(a,b;c;z),
\end{eqnarray*}
so that
$$
\frac{u(z)-v(z)}{z} = \frac{c-b}{c} F(a,b;c+1;z) -F(a,b;c;z) \to -b/c,$$
 as $z \to 0.$ Also, by (\ref{eq:viides}), 
$zv_1(z) \to 0$ as $z \to 0.$ Hence, letting $z \to 0,$ we get
\begin{eqnarray*}
A &=& (1-a)u_1(1) = (1-a)\frac{\Gamma(c)\Gamma(c+1-a-b)}{\Gamma(c+1-a)\Gamma(c-b)} \\
&= & (1-a)\frac{(\Gamma(c))^2}{ (1-a)\Gamma(a)\Gamma(b)} = \frac{(\Gamma(c))^2}{\Gamma(a)\Gamma(b)},
\end{eqnarray*}
 as claimed.

Note that, in Case (3), $ \Gamma(c) = \Gamma(1) = 1, \Gamma(b) = \Gamma(1-a),$
and thus by \cite[6.1.17]{AS}
$A= 1/(\Gamma(a)\Gamma(1-a)) = (\sin \pi a)/\pi .$ 
\ \ $\square$

%%%%%%%%%%%%%%%%%
%%%%%%%%%%%%%%%%%
%%%%%%%%%%%%%%%%%
\bigskip

For rational triples $(a,b,c)$ there are numerous cases
where the hypergeometric function $F(a,b;c;z)$ reduces to a simpler
function (see \cite{PBM}). Other important particular cases are
\emph{generalized elliptic integrals}, which we will now discuss.
For $a,r\in(0,1)$, the \emph{generalized elliptic integral
of the first kind} is given by
\begin{align*}
\mathcal{K}_a &= \mathcal{K}_a(r) = \frac{\pi}{2}
F(a,1-a;1;r^2)\\
&= (\sin\pi a)\int_0^{\pi /2}(\tan t)^{1-2a}(1-r^2\sin^2t)^{-a}\,
dt,\\
\mathcal{K}_a' &=\mathcal{K}_a'(r)=\mathcal{K}_a(r').
\end{align*}
We also define
\begin{equation*}
\mu_a(r)=\frac{\pi}{2\sin (\pi a)}
\frac{\mathcal{K}_a'(r)}{\mathcal{K}_a(r)},\quad r'=\sqrt{1-r^2}.
\end{equation*}

The {\it invariant} of the linear differential equation
\begin{equation}
\label{eq:wpw}
w''+pw'+qw=0,  
\end{equation}
where $p$ and $q$ are functions of $z$, is defined to be
$$
I\equiv q- \frac{1}{2}p'- \frac{1}{4}p^2
$$
(cf. [R2,p.9]).  If  $w_1$ and $w_2$ are two linearly independent 
solutions of (\ref{eq:wpw}), then their quotient 
$w\equiv w_2/w_1$ satisfies the differential equation
$$
S_w(z)=2I,
$$
where $S_w$ is the Schwarzian derivative
$$
S_w\equiv \left(\frac{w''}{w'}\right)'- 
\frac{1}{2}\left(\frac{w''}{w'}\right)^2
$$
and the primes indicate differentiations (cf. \cite[pp. 18,19]{R2}).  

From these considerations and the fact that $\mathcal {K}_a(r)$ 
and $\mathcal {K}'_a(r)$ are linearly independent solutions 
of (\ref {eq:kymm2}) (see \cite[(1.11)]{AQVV}), it follows that 
$w = \mu_a(r)$ satisfies the differential equation
$$
S_w(r) = \frac{-8a(1-a)}{(r')^2} + \frac{1+6r^2-3r^4}{2r^2(r')^4}.
$$

%%%%%%%%%%%%%%%%%%%%%%%%%%%%%%%%%%%

%%%%%%%%%%%%%%%%%%%%%%%%%%%
The \emph{generalized elliptic integral
of the second kind} is given by
\begin{align*}
\mathcal{E}_a &=\mathcal{E}_a(r)\equiv \frac{\pi}{2}
F(a-1,1-a;1;r^2)\\
&=(\sin \pi a)\int_0^{\pi/2}(\tan t)^{1-2a}(1-r^2\sin^2t)^{1-a}\,dt\\
\mathcal{E}_a' &=\mathcal{E}_a'(r)=\mathcal{E}_a(r'),\\
\mathcal{E}_a(0) &=\frac{\pi}{2},\quad \mathcal{E}_a(1)=
\frac{\sin (\pi a)}{2(1-a)}.
\end{align*}
For $a = \frac{1}{2}, ~\mathcal{K}_a$ and $\mathcal{E}_a$ 
reduce to $\mathcal{K}$ 
and $\mathcal{E}$, respectively, the usual elliptic
integrals of the first and second kind, respectively.
Likewise $\mu_{1/2}(r) = \mu(r)$, the modulus of the well-known 
Gr\"otzsch ring in the plane \cite{LV}.

\begin{cor}\label{cor:KaEa}
        The generalized elliptic integrals  $\mathcal {K}_a$ and 
$\mathcal {E}_a$ satisfy the differential equations
\begin{align}
&r(r')^2\frac{d^2\mathcal{K}_a}{dr^2} + (1-3r^2)
\frac{d\mathcal{K}_a}{dr}-4a(1-a)r\mathcal{K}_a=0,\label{eq:dereka}\\
&r(r')^2\frac{d^2\mathcal{E}_a}{dr^2}+(r')^2\frac{d\mathcal{E}_a}{dr}
+4(1-a)^2r\mathcal{E}_a=0,\label{eq:dertoka}
\end{align}
respectively.
\end{cor}

        {\bf Proof.}  These follow from (\ref{eq:kymm2}).\ \ $\square$
\medskip

For $a = \frac{1}{2}$ these reduce to well-known
differential equations \cite[pp.\ 474-475]{AVV5}, \cite{BF}.

\medskip

%\section{Hypergeometric Function and Generalized elliptic 
% integrals}\label{sect:genell}

\section{Identities of Legendre and Elliott}\label{sect:idenLE}

In geometric function theory the complete elliptic integrals 
$\mathcal{K}(r)$ and $\mathcal{E}(r)$ play an 
important role. These integrals may be defined, respectively, as
\begin{equation*}
\mathcal{K}(r)=
\textstyle{\frac{\pi}{2}} F(\textstyle{\frac{1}{2}},\textstyle{\frac{1}{2}}; 
1; r^2),\
\mathcal{E}(r)=
\textstyle{\frac{\pi}{2}} F(\textstyle{\frac{1}{2}},
-\textstyle{\frac{1}{2}};1;r^2),
\end{equation*}
for $-1 < r < 1$. These are $\mathcal{K}_a(r)$ and $\mathcal{E}_a(r)$, 
respectively, with $a=\frac{1}{2}$. We also consider the functions
\begin{align*}
\mathcal{K}' &=\mathcal{K}'(r)=\mathcal{K}(r'),\quad 0 < r < 1,\\
\mathcal{K}(0) &=\pi /2,\quad \mathcal{K}(1^-)=+\infty,
\end{align*}
and
\begin{equation*}
\mathcal{E}'=\mathcal{E}'(r)=\mathcal{E}(r'),\quad 0 \leqslant r 
\leqslant 1,
\end{equation*}
where $r'=\sqrt{1-r^2}$.
For example, these functions occur in the following quasiconformal counterpart
of the Schwarz Lemma \cite{LV}:

\begin{thm}\label{thm:kin}
For $K\in [1,\infty )$, let  $w$  be a $K$-quasiconformal mapping
of the unit disk $D=\{z:|z|<1\}$ into the unit disk $D'=\{w:|w|<1\}$ 
with $w(0)=0$.  Then
\begin{equation*}
|w(z)| \leqslant \varphi_K(|z|),
\end{equation*}
where
\begin{equation}\label{eqn:vphi}
\varphi_K(r)\equiv \mu^{-1}\big(\frac{1}{K}\mu(r)\big)\quad \text{and}
\quad \mu(r) \equiv \frac{\pi \mathcal{K}'(r)}{2 \mathcal{K}(r)}.
\end{equation}
This result is sharp in the sense that for each $z\in D$ and 
$K \in [1,\infty)$ 
there is an extremal $K$-quasiconformal mapping that takes the unit 
disk $D$ onto the unit disk $D'$ with $w(0)=0$ and 
$|w(z)|=\varphi_K(|z|)$ {\rm (}see \cite[p. 63]{LV}{\rm )}.
\end{thm}

It is well known \cite{BF} that the complete elliptic integrals 
$\mathcal{K}$ and $\mathcal{E}$ satisfy the Legendre relation
\begin{equation}\label{eq:legendre}
\mathcal{E} \mathcal{K}' + \mathcal{E}' \mathcal{K} - \mathcal{K} 
\mathcal{K}' = \frac{\pi}{2}.
\end{equation}
For several proofs of \eqref{eq:legendre} 
see \cite{Du}.

In 1904, E. B. Elliott \cite{E} (cf. \cite{AVV3}) obtained the following
generalization of this result.

\begin{thm}\label{thm:Ell} If  $a,b,c \geqslant 0$ 
and $0 < x < 1$ then
\begin{equation}\label{eq:elliott}
F_1F_2 + F_3F_4 - F_2F_3 = 
\frac{\Gamma(a+b+1)\Gamma(b+c+1)}{\Gamma (a+b+c+\frac{3}{2}) 
\Gamma(b + \frac{1}{2})}.
\end{equation}
where
\begin{align*}
F_1 &= F\bigg({\frac{1}{2}} + a, -\frac{1}{2} - c; 1+a+b; x\bigg),\\
F_2 &= F\bigg(\frac{1}{2} - a, \frac{1}{2} + c; 1+b+c; 1 - x\bigg),\\
F_3 &= F\bigg(\frac{1}{2} + a, \frac{1}{2} - c; 1+a+b;x\bigg),\\
F_4 &= F\bigg(-\frac{1}{2} - a, \frac{1}{2} + c; 1+b + c;1 - x\bigg).
\end{align*}
\end{thm}

Clearly \eqref{eq:legendre} is a special case of \eqref{eq:elliott},  
when $a = b = c = 0$ and $x=r^2$. 
%Duren \cite{Du} generalized \eqref{eq:legendre} by showing that if 
%$2c = a+b+1$ then both $F(x)$ and $F(1-x)$ are 
%solutions of (\ref{eq:hde}) as was mentioned in Section 4.  
%Other 
For a discussion of generalizations of Legendre's 
Relation see Karatsuba and Vuorinen 
\cite{KV} and Balasubramanian, Ponnusamy, Sunanda Naik, and Vuorinen
\cite{BPSV}.

Elliott proved \eqref{eq:elliott} by a clever change of variables in 
multiple integrals.  Another proof was suggested without details in 
\cite[p.\ 138]{AAR}, and here we provide the missing details.

{\bf Proof of Theorem \ref{thm:Ell}.}  In particular, let $y_1 \equiv F_3$, 
$y_2 \equiv x^{-a-b}(1-x)^{b+c} F_2$. Then by \cite[pp.\ 174, 175]{R2} or 
\cite[(3.2.12), (3.2.13)]{AAR}, $y_1$ and $y_2$ are linearly 
independent solutions of \eqref{eq:kymm}.

  By \eqref{eq:yhd},
\begin{equation}\label{eq:by38}
x (1\!-\!x)y'_1 = \big(a+b+c+\frac{1}{2}\big) F_1 +
[-(a+b+c+\frac{1}{2}) + (a + \frac{1}{2}) x] F_3,
\end{equation}
and by  \eqref{eq:seits},
\begin{equation}\label{eq:bythm}
\begin{split}
x(1-x) y'_2 = &x(1-x) \big[-(a+b) x^{-a-b-1} (1-x)^{b+c}\\
&-(b+c) x^{-a-b} (1-x)^{b+c-1} \big] F_2\\
&- x^{-a-b}(1-x)^{b+c} \big[\big(a+b+c+\frac{1}{2}\big) F_4\\
&+ \big[ - (a+b+c+\frac{1}{2} \big) + \big( c + \frac{1}{2} 
\big) (1-x)) F_2 \big].
\end{split}
\end{equation}
Multiplying \eqref{eq:bythm} by $y_1$ and \eqref{eq:by38} by $y_2$ 
and subtracting, we obtain
\begin{align*}
&x(1-x)(y_2 y'_1 - y_1 y'_2) = \big(a+b+c+\frac{1}{2}\big) x^{-a-b} 
(1-x)^{b+c} F_1 F_2 \\
&\quad+ \big[-\big(a+b+c+\frac{1}{2} \big)+ \big(a + \frac{1}{2} \big) x 
\big] x^{-a-b} (1-x)^{b+c} F_2 F_3 \\
&\quad+ x(1-x) [(a+b)x^{-a-b-1}(1-x)^{b+c}\\
&\quad+ (b+c) x^{-a-b} (1-x)^{b+c-1}] F_2F_3\\
&\quad+ x^{-a-b} (1-x)^{b+c} \big(a+b+c+\frac{1}{2}\big) F_3 F_4 \\
&\quad+ x^{-a-b}(1-x)^{b+c} \big[- (a+b+c+\frac{1}{2}\big) + \big(c + 
\frac{1}{2} \big) (1-x)\big] F_2 F_3\\
&= \big(a+b+c+\frac{1}{2}\big) x^{-a-b} (1-x)^{b+c} F_1 F_2\\
&\quad+ x^{-a-b}(1-x)^{b+c} \big[ -(a+b+c+\frac{1}{2}\big)\\ 
&\quad+ \big(a + \frac{1}{2}\big) x + (a+b) (1-x) \\
&\quad+ (b+c) x -\big( a+b+c+\frac{1}{2}\big) + \big(c + \frac{1}{2}\big) 
(1-x) \big] F_2 F_3 \\
&\quad+ x^{-a-b} (1-x)^{b+c} \big(a+b+c+ \frac{1}{2}\big) F_3 F_4 \\
&= \big(a+b+c+\frac{1}{2} \big) x^{-a-b} (1-x)^{b+c} F_1 F_2 \\
&\quad+ x^{-a-b} (1-x)^{b+c} \big[ - \big(a+b+c+\textstyle{\frac{1}{2}}
\big) \big]F_2 F_3 \\
&\quad+ \big(a+b+c+\frac{1}{2}\big) x^{-a-b} (1-x)^{b+c} F_3F_4.
\end{align*}
So
\begin{align*}
x(1-x)&W(y_2,y_1) = \frac{A}{x^{a+b} (1-x)^{-b-c}}\\
&= x^{-a-b} (1-x)^{b+c}\big( a+b+c+\frac{1}{2}\big) 
[F_1 F_2 + F_3 F_4 - F_2 F_3]
\end{align*}
by Lemma \ref{lem:linind}. Thus
\begin{equation*}
F_1 F_2 + F_3 F_4 - F_2 F_3 = A,
\end{equation*}
where $A$ is a constant.

\medskip

Now, by \eqref{eq:viides},
\begin{equation*}
F_1 F_2 \mbox{ tends to } F \big(\frac{1}{2} + a, - \frac{1}{2} - c; 
a+b+1;1\big) =
\frac{\Gamma(a+b+1)\Gamma (b+c+1)}{\Gamma(b+ \frac{1}{2})
\Gamma(a+b+c+\frac{3}{2})}
\end{equation*}
as $x \rightarrow 1$, since $\frac{1}{2} + a + (-\frac{1}{2} - c) 
= a - c < a +b+1$.

Next
\begin{equation*}
F_3 F_4 - F_3 F_2 = F_3(F_4 - F_2),
\end{equation*}
where $F_4 - F_2 \sim\, \text{const}\,\cdot\,(1-x)^2 + O((1-x)^3)$, 
and
\begin{equation*}
F_3 = \frac{\Gamma(a+b+1) \Gamma(b+c+1)}{\Gamma (b+ \frac{1}{2}) 
\Gamma (a+b+c+\frac{1}{2})}
\end{equation*}
if $a + \frac{1}{2} + \frac{1}{2} - c < a + b + 1$, or $-c < b$, 
i.e., $b > 0$ or $c>0$.  If $-c = b = 0$, then
\begin{equation*}
F_3 = \frac{R(a + \frac{1}{2},\frac{1}{2}-c)}{B(a + \frac{1}{2}, 
\frac{1}{2}-c)} + O((1-x) \log(1-x))
\end{equation*}
by \eqref{eq:viides}.
%Here $\gamma$ is the Euler-Mascheroni constant and
%$\Psi(x) \equiv \Gamma'(x)/\Gamma(x)$.
In either case the product $F_3(F_4 - F_2)$
tends to $0$ as $x \rightarrow 1$.
The third case $a + \frac{1}{2} + \frac{1}{2} - c> a + b + 1$
is impossible since we are assuming that $b,c$ are nonnegative.
Thus $A = \Gamma (a+b+1) \Gamma (b+c+1)/(\Gamma (b + \frac{1}{2})
\Gamma (a+b+c+\frac{3}{2}))$, as desired.\ \ $\square$

\bigskip

The generalized elliptic integrals satisfy the identity
\begin{equation}
\label{eq:legkaea}
\mathcal{E}_a\mathcal{K}_a'+\mathcal{E}_a'\mathcal{K}_a-
\mathcal{K}_a\mathcal{ K}_a'=\frac{\pi\sin(\pi a)}{4(1-a)}.
\end{equation}
This follows from Elliott's formula \eqref{eq:elliott} and  
contains the classical relation of Legendre \eqref{eq:legendre} 
as a special case.

Finally, we record the following formula of
Kummer \cite[p. 63, Form. 30]{Kum}:
$$
\hspace{-1in}F(a,b;a+b-c+1;1-x)F(a+1,b+1;c+1;x) 
$$
$$
\qquad \qquad + \frac{c}{a+b-c+1}F(a,b;c;x)F(a+1,b+1;a+b-c+2;1-x)
$$
$$
\hspace{.25in}  = Dx^{-c}(1-x)^{c-a-b-1},\ \ D= \frac{\Gamma(a+b-c+1) \Gamma(c+1)}{\Gamma(a+1) \Gamma(b+1)}.
$$
This formula, like Elliott's identity, may be rewritten in many
different ways if we use the contiguous relations of Gauss.
Note also the special case $c=a+b-c+1.$

\bigskip

\section{Mean values}\label{sect:meanval}

The \emph{arithmetic-geometric mean} of positive numbers $a, b>0$ is
the limit
\begin{equation*}
AGM(a,b)= \lim a_n = \lim b_n,
\end{equation*}
%where $a_0=a$, $b_0=b$, $a_{n+1}=(a_n+b_n)/2$, 
%$b_{n+1}=\sqrt{a_n b_n}$. 
 where $a_0=a$, $b_0=b$, and for $n=0,1,2,3,...,$
 $$a_{n+1}=A(a_n,b_n) \equiv (a_n+b_n)/2, 
 \quad b_{n+1}=G(a_n, b_n)\equiv\sqrt{a_n b_n},$$ 
 are the arithmetic and geometric means of $a_n$ and $ b_n,$ resp.
For a mean value $M$, we also 
consider the $t$-\emph{modification} defined as
\begin{equation*}
M_t(a,b)= M(a^t, b^t)^{1/t}.
\end{equation*}
For example, the \emph{power mean} of $a,b >0$ is
\begin{equation*}
A_t(a,b) =\left( \frac{a^t +b^t}{2}\right)^{1/t},
\end{equation*}
and the \emph{logarithmic mean} is
\begin{equation*}
L(a,b) = \frac{a-b}{\log (b/a)}.
\end{equation*}
The power mean is the $t$-modification of the arithmetic mean  $A_1(a,b)$.

\medskip

The connection between mean values and elliptic integrals is
provided by Gauss's amazing result
\begin{equation*}
AGM(1, r') = \frac{\pi}{2\mathcal{K}(r)}.\end{equation*}
This formula motivates  the question of finding minorant/majorant 
functions for $\mathcal{K}(r)$ in terms of mean values.
For a fixed $x>0$ the function $t  \mapsto L_t(1,x), t >0,$ increases
with $t$ by \cite[Theorem 1.2 (1)]{VV}. The two-sided inequality
$$ L_{3/2}(1,x) > AG(1,x)> L(1,x) $$
holds; the second inequality was pointed out in \cite{CV}, and the first
one, due to J. and P. Borwein \cite{BB2}, proves a sharp 
estimate settling a question raised in connection with \cite{VV}.
Combined with the identity above, this inequality yields
a very precise inequality for $\mathcal{K}(r).$

Several inequalities between mean values have been proved
recently. See, for instance, \cite{AlQ2}, \cite{QS}, \cite{S1}, 
\cite{S2}, \cite{S3}, \cite{T}, \cite{C},  and  \cite{Br}.

Finally, we remark that the arithmetic-geometric mean, together with
Legendre's Relation, played a central role in a rapidly
converging algorithm for the number $\pi$ in \cite{Sa}.
See also \cite{BB1, H, Le, Lu}.

\medskip

\section{Approximation of elliptic integrals}\label{sect:apprell}

Efficient algorithms for the numerical
evaluation of $\mathcal{K}(r)$ and $\mathcal{E}(r)$
are based on the arithmetic-geometric mean iteration of Gauss.
This fact led to some
close majorant/minorant functions for $\mathcal{K}(r)$ in terms
of mean values in \cite{VV}. 

Next, let $a$ and $b$ be the semiaxes of an ellipse with 
$a >b$ and eccentricity 
$e = \sqrt{a^2 - b^2}/a$, and let $L(a,b)$ denote the arc length of 
the ellipse. Without loss of generality we take $a = 1$.  
In 1742, Maclaurin (cf. \cite{AB}) determined that
\begin{equation*}
L(1,b)=4\mathcal{E}(e)=2\pi\cdot {}_2F_1(\textstyle{\frac{1}{2}},
-\textstyle{\frac{1}{2}};1;e^2).
\end{equation*}

In 1883, Muir (cf. \cite{AB}) proposed that $L(1,b)$ could be 
approximated by the expression $2\pi [(1+b^{3/2})/2]^{2/3}$.
Since this expression has a close resemblance to the power
mean values studied in \cite{VV}, it is natural to study
the sharpness of this approximation.
Close numerical examination of the error in this
approximation led Vuorinen \cite{V2} to conjecture that Muir's
approximation is a lower bound for the arc length.  Letting
$r=\sqrt{1-b^2}$, Vuorinen asked whether
\begin{equation}\label{eq:vuor}
\frac{2}{\pi}\mathcal{E}(r) 
= {}_2F_1 \Big(\tfrac{1}{2},- \tfrac{1}{2}; 1;r^2\Big)
\geqslant \Big( \frac{1+(r')^{3/2}}{2} \Big)^{2/3}
\end{equation}
for all $r\in [0,1]$.

In \cite{BPR1} Barnard and his coauthors proved that
inequality \eqref{eq:vuor} is true.
In fact, they expanded both functions into Maclaurin
series and proved that the differences of the corresponding 
coefficients of the two series all have the same sign.

Later, the same authors \cite{BPR2} discovered an upper bound for
$\mathcal{E}$ that complements the lower bound in \eqref{eq:vuor}:
\begin{equation}\label{eq:upper}
\frac{2}{\pi}\mathcal{E}(r)
= {}_2F_1 \Big(\tfrac{1}{2}, - \tfrac{1}{2}; 1;r^2\Big) 
\leqslant \Big( \frac{1 +(r')^2}{2} \Big)^{1/2}, 
\quad 0 \leqslant r \leqslant 1.
\end{equation}
See also \cite{BPS}.

In \cite{BPR2} the authors have considered 13 historical 
approximations (by Kepler, Euler, Peano, Muir, Ramanujan, 
and others) for the arc length of an ellipse and determined a 
linear ordering among them. Their main tool was the following 
Lemma \ref{lem:genhyper} on generalized hypergeometric functions.
These functions are
defined by the formula
\begin{equation*}
_pF_q(a_1,a_2,\cdots ,a_p;b_1,b_2,\cdots ,b_q;z)\equiv 
1+\sum_{n=1}^{\infty}\frac{\Pi_{i=1}^p(a_i,n)}{\Pi_{j=1}^q(b_j,n)}
\cdot \frac{z^n}{n!},
\end{equation*}
where  $p$ and $q$ are positive integers and in which no 
denominator parameter $b_j$ is permitted to
be zero or a negative integer. When $p=2$ and $q=1$, this reduces to 
the usual Gaussian hypergeometric function $F(a,b;c;z)$.

\begin{lem}\label{lem:genhyper}
Suppose $a,b > 0$. Then for any $\epsilon$ satisfying
$\frac{ab}{1 + a + b} < \epsilon < 1$,
\begin{equation*}
_3F_2(-n,a,b; 1 + a + b, 1 + \epsilon - n;1) > 0
\end{equation*}
for all integers  $n \geqslant 1$.
\end{lem}

\begin{OtherAppr}
At the end of the preceding section we pointed out that
upper and lower bounds can be found for $\mathcal{K}(r) $
in terms of mean values.
Another source for the approximation of $\mathcal{K}(r)$ is 
based on the asymptotic behavior at the singularity $r=1$, 
where $\mathcal{K}(r)$ has logarithmic growth. Some of the 
approximations motivated by this aspect will be discussed next.

Anderson, Vamanamurthy, and Vuorinen  \cite{AVV2} approximated 
$\mathcal{K}(r)$ by the inverse hyperbolic tangent function 
$\arth$, obtaining the inequalities
\begin{equation}\label{eq:arth}
\frac{\pi}{2}\Bigg(\frac{\arth r}{r}\Bigg)^{1/2} 
< \mathcal{K}(r) 
< \frac{\pi}{2}\, \frac{\arth r}{r},
\end{equation}
for $0<r<1$. Further results were proved by Laforgia and Sismondi 
\cite{LS}. K\"uhnau \cite{Ku} and Qiu \cite{Q1} proved that, for 
$0 < r < 1$, 
\begin{equation*}
\frac{9}{8+r^2} < \frac{\mathcal{K}(r)}{\log (4/r')}.
\end{equation*}

Qiu and Vamanamurthy \cite{QVa} proved that
\begin{equation*}
\frac{\mathcal{K}(r)}{\log (4/r')} 
< 1 + \frac{1}{4}(r')^2\quad \text{for}\ 0 < r < 1.
\end{equation*}
Several inequalities for $\mathcal{K}(r)$ are given in 
\cite[Theorem 3.21]{AVV5}. 
Later Alzer \cite{A3} showed that
\begin{equation*}
1+\Big(\frac{\pi}{4\log 2}-1\Big)(r')^2
<\frac{\mathcal{K}(r)}{\log (4/r')}, 
\end{equation*}
for $0<r<1$. He also showed that the constants
$\frac{1}{4}$ and $\pi/( 4\log 2)-1$ in the above 
inequalities are best possible.

For further refinements, see \cite[(2.24)]{QVu1} and \cite{Be}.

Alzer and Qiu \cite{AlQ1} have written a related manuscript in 
which, besides proving many inequalities for complete elliptic 
integrals, they have refined \eqref{eq:arth} by proving that
\begin{equation*}
\frac{\pi}{2}\Big(\frac{\arth r}{r}\Big)^{3/4}<\mathcal{K}(r) 
< \frac{\pi}{2}\, \frac{\arth r}{r}.
\end{equation*}
They also showed that $3/4$ and $1$ are the best exponents for 
$(\arth r)/r$ on the left and right, respectively.

One of the interesting tools of these  authors is the following lemma of 
Biernaki and Krzy\.{z} \cite{BK} (for a detailed proof see 
\cite{PV1}):
\end{OtherAppr}

\begin{lem}\label{lem:BK}
Let $r_n$ and $s_n$, $n = 1,2,\dots$ be real numbers, and let
the power series $R(x) = \sum_{n=1}^{\infty }r_nx^n$ and $S(x) =
\sum_{n=1}^{\infty}s_nx^n$ be convergent for
$|x|<1$. If $s_n > 0$ for $n=1,2,\ldots$, and if
$r_n/s_n$ is strictly increasing {\rm(}resp. decreasing{\rm)} for 
$n=1,2,\ldots$, then the function $R/S$ is strictly increasing 
{\rm(}resp. decreasing{\rm)} on $(0,1)$.
\end{lem}

\begin{GenEllRmk} 
For the case of generalized elliptic
integrals some inequalities are given in \cite{AQVV}.
B. C. Carlson has introduced some standard forms for elliptic 
integrals involving certain symmetric
integrals. Approximations for these functions can be found in 
\cite{CG}. 
\end{GenEllRmk}

\section{Landen inequalities}

It is well known (cf. \cite{BF}) that the complete elliptic 
integral of the first kind satisfies the Landen identities
\begin{equation*}
\mathcal{K}\left(\frac{2\sqrt r}{1+r}\right)=(1+r)\mathcal{K}(r),
\quad \mathcal{K}\left(\frac{1-r}{1+r}\right)
=\frac{1+r}{2}\mathcal{K}'(r).
\end{equation*}
Recall that $\mathcal{K}(r)=\frac{\pi}{2}F(\frac{1}{2}$, 
$\frac{1}{2};1;r^2)$.
It is thus natural to consider, as suggested in \cite{AVV3},
the problem of finding an analogue
of these formulas for the zero-balanced hypergeometric
function $F(a,b;c;r)$ for $a,b,c>0$ and $a+b=c$, at least when
the parameters $(a,b,c)$ are close to $(\frac{1}{2}$, 
$\frac{1}{2},1)$. From (3.3) it is clear that $F(a, b;c;r^2)$
has a logarithmic singularity at $r=1$,  if $a,b>0$, 
$c=a+b$ (cf. \cite{AAR}). Some refinements of the growth estimates were given 
in \cite{ABRVV} and \cite{PV1}.

Qiu and Vuorinen \cite{QVu1} proved the following  Landen-type 
inequalities: For $a,b\in (0,1)$, $c=a+b$,
\begin{align*}
F\bigg(a,b;c;\bigg(\frac{2\sqrt r}{1+r}\bigg)^2\bigg)
&\leqslant (1+r)F(a,b;c;r^2)\\
&\leqslant F\bigg(a,b;c;\bigg(\frac{2\sqrt r}{1+r}\bigg)^2\bigg)
+\frac{1}{B}(R-\log 16)
\end{align*}
and
\begin{align*}
\frac{1+r}{2}F(a,b;c;1-r^2)
&\leqslant F\bigg(a,b;c;\bigg(\frac{1-r}{1+r}\bigg)^2\bigg)\\
&\leqslant \frac{1+r}{2}\bigg[F(a,b;c;1-r^2)+\frac{1}{B}
(R-\log 16)\bigg],
\end{align*}
with equality in each instance if and only if $a=b=\frac{1}{2}$.  
Here $B=B(a,b)$, the beta function, and 
$R=R(a,b)=-2\gamma-\Psi (a)-\Psi (b)$, where $\Psi$ is as given in Section 2.
%Here $B$ denotes the beta function $B(a,b)\equiv \Gamma (a)\Gamma
%(b)/\Gamma (a+b)$ and $R=R(a,b)\equiv -\Psi (a)-\Psi (b)-2\gamma$, where
%$\Psi =\Gamma'/\Gamma$ and $\gamma$ is the Euler-Mascheroni constant.

\section{Hypergeometric series as an analytic function}

For rational triples $(a,b,c)$ the hypergeometric function
often can be expressed in terms of elementary functions.
Long lists with such triples containing hundreds of functions
can be found in \cite{PBM}. For example, the functions
$$ f(z)\equiv z F(1,1;2;z)= -\log(1-z)$$
and
$$g(z)\equiv z F(1,1/2;3/2;z^2)=
 \frac{1}{2} \log \left(\frac{1+z}{1-z}\right)$$
have the property that they both map the unit disk into a strip domain.
Observing that they both correspond to the case $c=a+b$ one may ask
(see \cite{PV1,PV2}) whether there exists $\delta >0$ such that
$zF(a,b;a+b;z)$ and $zF(a,b;a+b;z^2)$ with $a,b \in (0,\delta)$
map into a strip domain.

Membership of hypergeometric functions in some special classes of
univalent functions is studied in \cite{PV1,PV2,BPV2}.

\vskip 1cm

\noindent
ANDERSON: \\
  Department of Mathematics \\
Michigan State University \\
     East Lansing, MI 48824, USA \\
      email: {\tt anderson@math.msu.edu}\\
     FAX: +1-517-432-1562\\ [1mm]

%\noindent
%QIU:\\
%    Hangzhou Institute of Electronics Engineering\\
%     Hangzhou 310037, P. R. CHINA\\
%     FAX: 086-571-8077232\\[1mm]

\noindent
VAMANAMURTHY:\\
   Department of Mathematics \\
    University of Auckland \\
    Auckland, NEW ZEALAND\\
     email: {\tt vamanamu@math.auckland.nz}\\
FAX: +649-373-7457\\ [2mm]

\noindent
VUORINEN:\\
     Department of Mathematics \\
     University of Helsinki \\
     P.O. Box 4 (Yliopistonkatu 5)\\
     FIN-00014, FINLAND\\
     e-mail: ~~{\tt vuorinen@csc.fi}\\
     FAX: +358-9-19123213\\


\begin{thebibliography}{AAAAA}

\bibitem[AS]{AS} \textsc{M. Abramowitz and I. A. Stegun, eds:}
\emph{Handbook of Mathematical Functions with Formulas, Graphs and
Mathematical Tables}, Dover, New York, 1965.

\bibitem[AB]{AB} \textsc{G. Almkvist and B. Berndt:} Gauss, Landen, Ramanujan,
the arithmetic-geometric mean, ellipses, pi, and the
Ladies Diary, \emph{Amer. Math. Monthly} \textbf{95} (1988), 585--608.

\bibitem[A1]{A1} \textsc{H. Alzer:} On some inequalities for 
the gamma and psi 
functions, \emph{Math. Comp.} \textbf{66} (1997), 373--389.

\bibitem[A2]{A2} \textsc{H. Alzer:} Inequalities for the gamma 
and polygamma functions, \emph{Abh. Math. Sem. Univ. 
Hamburg} \textbf{68}, (1998), 363--372.

\bibitem[A3]{A3} \textsc{H. Alzer:} 
Sharp inequalities for the complete elliptic 
integral of the first kind, \emph{Math. Proc. Camb. Phil. Soc.} 
\textbf{124} (1998), 309--314.

\bibitem[A4]{A4} \textsc{H. Alzer:} Inequalities for the gamma function, 
\emph{Proc. Amer. Math. Soc.} \textbf{128} (2000), 141--147.

\bibitem[A5]{A5} \textsc{H. Alzer:} Inequalities for the volume of the unit 
ball in $\R^n$, \emph{J. Math. Anal. Appl.} \textbf{252} (2000), 353--363.

\bibitem[A6]{A6} \textsc{H. Alzer:} A power mean inequality for the gamma 
function, \emph{Monatsh. Math.} \textbf{131} (2000), 179--188.

\bibitem[A7]{A7} \textsc{H. Alzer:} Inequalities for the Hurwitz zeta 
function, \emph{Proc. Royal Soc. Edinb.} \textbf{130A} (2000), 1227--1236.

\bibitem[A8]{A8} \textsc{H. Alzer:} Mean value inequalities for the polygamma 
function, \emph{Aeqationes Math.} \textbf{61} (2001), 151--161.

\bibitem[A9]{A9} \textsc{H. Alzer:} {Sharp inequalities for the beta 
function}, \emph{ Indag. Math. (N.S.)} \textbf{12} (2001), 15--21.

\bibitem[AlQ1]{AlQ1} \textsc{H. Alzer and S.-L. Qiu:} 
{Monotonicity theorems and inequalities for the complete elliptic integrals}, 
Manuscript (2000).

\bibitem[AlQ2]{AlQ2} \textsc{H. Alzer and  S.-L. Qiu:} 
{Inequalities for means in two variables}, Manuscript (2000).




\bibitem[AW]{AW} \textsc{H. Alzer and J. Wells:} Inequalities for 
the polygamma functions, \emph{SIAM J. Math. Anal.} \textbf{29} 
(1998), 1459--1466 (electronic).

\bibitem[ABRVV]{ABRVV} \textsc{G. D. Anderson, R. W. Barnard, K. C. Richards, 
M. K. Vamanamurthy, and M. Vuorinen:} Inequalities for zero-balanced 
hypergeometric functions, \emph{Trans. Amer. Math. Soc.} \textbf{347} 
(1995), 1713--1723.

\bibitem[AnQ]{AnQ} \textsc{G. D. Anderson and S.-L. Qiu:} 
A monotoneity property of the gamma function, 
\emph{Proc. Amer. Math. Soc.} \textbf{125} (1997), 3355--3362.

\bibitem[AQVV]{AQVV} \textsc{G. D. Anderson, S.-L. Qiu, M. K. Vamanamurthy, 
and M. Vuorinen:} Generalized elliptic integrals and modular equations, 
\emph{Pacific J. Math.} \textbf{192} (2000), 1--37.


\bibitem[AQVa]{AQVa} \textsc{G. D. Anderson, S.-L. Qiu, M. K. Vamanamurthy:} 
Elliptic integral inequalities, with applications. 
\emph{Constr. Approx.}  \textbf{14} (1998), no. 2, 195--207. 


\bibitem[AQVu]{AQVu} \textsc{G. D. Anderson, S.-L. Qiu, and M. Vuorinen:}
Precise estimates for differences of the Gaussian hypergeometric function,
\emph{J. Math. Anal. Appl.}  \textbf{215} (1997), 212--234.


\bibitem[AVV1]{AVV1} \textsc{G. D. Anderson, M. K. Vamanamurthy, and 
M. Vuorinen:} Special functions of quasiconformal theory, 
\emph{Exposition. Math.} \textbf{7} (1989), 97--138.

\bibitem[AVV2]{AVV2} \textsc{G. D. Anderson, M. K. Vamanamurthy, and 
M. Vuorinen:} Functional inequalities for hypergeometric functions and 
complete elliptic integrals, \emph{SIAM J. Math. Anal.} \textbf{23} 
(1992), 512--524.

\bibitem[AVV3]{AVV3} \textsc{G. D. Anderson, M. K. Vamanamurthy, and 
M. Vuorinen:} Hypergeometric functions and elliptic integrals, 
in \emph{Current Topics in Analytic Function Theory}, ed. by 
H. M. Srivastava and S. Owa, World Scientific, London, 1992, 
pp. 48--85.

\bibitem[AVV4]{AVV4} \textsc{G. D. Anderson, M. K. Vamanamurthy, and 
M. Vuorinen:} Inequalities for quasiconformal mappings in space,
\emph{Pacific J. Math.} \textbf{160} (1993), 1--18.

\bibitem[AVV5]{AVV5} \textsc{G. D. Anderson, M. K. Vamanamurthy, and 
M. Vuorinen:} \emph{Conformal Invariants, Inequalities, and 
Quasiconformal Maps}, J. Wiley, 1997.

\bibitem[AAR]{AAR} \textsc{G. Andrews, R. Askey, R. Roy:} \emph{Special 
Functions,} Encyclopedia of Mathematics and its Applications,
Vol. 71, Cambridge U. Press, 1999.

\bibitem[As]{As} \textsc{R. Askey:} Ramanujan and hypergeometric and 
basic hypergeometric series, Ramanujan Internat. Symposium on 
Analysis, December 26-28, 1987, ed. by N. K. Thakare, 1-83, 
Pune, India, \emph{Russian Math. Surveys} \textbf{451} (1990), 37-86.

\bibitem[BPV1]{BPV1} \textsc{R. Balasubramanian, S. Ponnusamy, and 
M. Vuorinen:} Functional inequalities for the quotients of 
hypergeometric functions, \emph{J. Math. Anal. Appl.} \textbf{218} 
(1998), 256--268.

\bibitem[BPV2]{BPV2} \textsc{R. Balasubramanian, S. Ponnusamy, and 
M. Vuorinen:} On hypergeometric functions and function spaces, 
\emph{J. Comp. Appl. Math.} (to appear).

\bibitem[BPSV]{BPSV}  \textsc{R. Balasubramanian, S. Ponnusamy, 
Sunanda Naik, and M. Vuorinen: }
Elliott's identity and hypergeometric functions,  Preprint 284,
April 2001, University of Helsinki, 22 pp.



\bibitem[B]{B} \textsc{R. W. Barnard:} On applications of hypergeometric 
functions, Continued fractions and geometric function theory 
(CONFUN)(Trondheim, 1997), \emph{J. Comp. Appl. Math.} \textbf{105}  
(1999), no. 1--2, 1--8.

\bibitem[BPR1]{BPR1} \textsc{R. W. Barnard, K. Pearce, and K. C. Richards:} 
A monotonicity property involving $_3F_2$ and comparisons of the 
classical approximations of elliptical arc length, 
\emph{SIAM J. Math. Anal.} \textbf{32} (2000), 403--419 (electronic).

\bibitem[BPR2]{BPR2} \textsc{R. W. Barnard, K. Pearce, and K. C. Richards:} 
An inequality involving the generalized
hypergeometric function and the arc length of an ellipse, 
\emph{SIAM J. Math. Anal.} \textbf{31} (2000), no. 3, 693--699 
(electronic).

\bibitem[BPS]{BPS} \textsc{R. W. Barnard, K. Pearce, and L. Schovanec:} 
Inequalities for the perimeter of an ellipse, 
\emph{J. Math. Anal. Appl.} \textbf{260} (2001), 295--306.


\bibitem[Be]{Be} \textsc{A. F. Beardon:} 
The hyperbolic metric of a rectangle,
\emph{ Ann. Acad. Sci. Fenn. Ser. A I} \textbf{26} (2001), 401--407. 


\bibitem[BP]{BP} \textsc{C. Berg and H. Pedersen:} A completely monotone
function related to the gamma function, \emph{J. Comp. Appl. Math.} 
\textbf{133} (2001), 219--230.

\bibitem[Be1]{Be1} \textsc{B. C. Berndt:} \emph{Ramanujan's Notebooks}, 
Vol. II, Springer-Verlag, Berlin, 1989.

\bibitem[Be2]{Be2} \textsc{B. C. Berndt:} \emph{Ramanujan's Notebooks}, 
Vol. IV, Springer-Verlag, Berlin, 1993.

\bibitem[BBG]{BBG} \textsc{B. C. Berndt, S. Bhargava, and F. G. Garvan:} 
Ramanujan's theories of elliptic functions to alternative bases, 
\emph{Trans. Amer. Math. Soc.} \textbf{347} (1995), 4163--4244.

\bibitem[BCK]{BCK} \textsc{B. C. Berndt, Y.-S. Choi, and S.-Y. Kang:}
The problems submitted by Ramanujan to the Journal
of the Indian Mathematical Society. Continued fractions: 
from analytic number theory to constructive approximation 
(Columbia, MO, 1998), 15--56, \emph{Contemp. Math.},  \textbf{236}, 
Amer. Math. Soc., Providence, RI, 1999.

\bibitem[BK]{BK} \textsc{M. Biernaki and J. Krzy\.z:} On the monotonicity 
of certain functionals in the theory of analytic functions,
\emph{Ann. Univ. M. Curie-Sk{\l}odowska} \textbf{2} (1955), 135--147.

\bibitem[BB1]{BB1} \textsc{J. M. Borwein and P. B. Borwein:} 
\emph{Pi and the AGM}, Wiley, New York, 1987.

\bibitem[BB2]{BB2} \textsc{J. M. Borwein and P. B. Borwein:} Inequalities 
for compound mean iterations with logarithmic asymptotes, 
\emph{J. Math. Anal. Appl.} \textbf{177} (1993), 572--582.

\bibitem[Br]{Br} \textsc{P. Bracken:} An arithmetic-geometric mean
inequality. \emph{Exposition. Math.}  \textbf{19} (2001), 273--279.

\bibitem[BF]{BF} \textsc{P. F. Byrd and M. D. Friedman:} \emph{Handbook of 
Elliptic Integrals for Engineers and Scientists}, 2nd ed., 
Grundlehren Math. Wiss., Vol. 67, Springer-Verlag, Berlin, 1971.

\bibitem[CG]{CG} \textsc{B. C. Carlson and J. L. Gustafson:}
Asymptotic approximations for symmetric elliptic integrals,
\emph{SIAM J. Math. Anal.} \textbf{25} (1994), 288--303.


\bibitem[CV]{CV} \textsc{B. C. Carlson and M. Vuorinen:} 
An inequality of the AGM and the logarithmic mean, 
\emph{SIAM Rev.} \textbf{33} (1991), Problem 91-17, 655.


\bibitem [Ch]{Ch}   \textsc{I. Chavel:} 
\emph{Riemannian Geometry -- A Modern Introduction,}
 Cambridge Tracts in Math. 108, Cambridge Univ. Press, 1993.


\bibitem[C]{C} \textsc{S.-Y. Chung:} Functional means and harmonic functional 
means, \emph{Bull. Austral. Math. Soc.} \textbf{57} (1998), 207--220.

\bibitem[De]{De} \textsc{D. W. DeTemple:} A quicker convergence to Euler's 
constant, \emph{Amer. Math. Monthly} \textbf{100(5)} (1993), 468--470.

\bibitem[Du]{Du} \textsc{P. L. Duren:} The Legendre relation for elliptic 
integrals, in \emph{Paul Halmos: Celebrating 50 years of Mathematics}, 
ed. by J. H. Ewing and F. W. Gehring, Springer-Verlag, New York, 
1991, pp. 305--315.

\bibitem[EL]{EL} \textsc{\'A. Elbert and A. Laforgia:} On some properties of 
the gamma function, \emph{Proc. Amer. Math. Soc.} \textbf{128} (2000), 
2667--2673.

\bibitem[E]{E} \textsc{E. B. Elliott:} A formula including Legendre's 
$EK' + KE' - KK' = \tfrac{1}{2}\pi$, \emph{Messenger of Math.} 
\textbf{33} (1904), 31--40.


\bibitem[G]{G} \textsc{W. Gautschi:}
The incomplete gamma functions since Tricomi. \emph{ Tricomi's ideas 
and contemporary applied mathematics (Rome/Turin, 1997),}
203--237, Atti Convegni Lincei, 147, Accad. Naz. Lincei, Rome, 1998. 

%\bibitem[HP]{HP} \textsc{J. Hersch and A. Pfluger:} Generalisation du lemme 
%de Schwarz et du principe de la mesure harmonique pour les fonctions 
%pseudo-analytiques, \emph{C. R. Acad. Sci. Paris} \textbf{234} (1952), 
%43--45.


\bibitem[H]{H} \textsc{M. Hirschhorn:} A new formula for $\pi,$
\emph{ Gaz. Austral. Math. Soc.} \textbf{25} (1998), 82--83.
%Manuscript 1997. %{\bf From mathematical constants www-page}



\bibitem[K1]{K1} \textsc{E. A. Karatsuba:} On the computation of the Euler 
constant $\gamma$, \emph{Numer. Algorithms} \textbf{24} (2000), 83--87.

\bibitem[K2]{K2} \textsc{E. A. Karatsuba:} {On the asymptotic 
representation of the Euler gamma function by Ramanujan}, 
 \emph{J. Comp. Appl. Math.} \textbf{135.2} (2001), 225--240.

%Preprint 248, University of Helsinki, 1999, 1--21.

\bibitem[KV]{KV} \textsc{E. A. Karatsuba and M. Vuorinen:} On
hypergeometric functions and generalizations of Legendre's relation,
\emph{J. Math. Anal. Appl.}  \textbf{260} (2001), 623-640.



\bibitem[KlR]{KlR} \textsc{D. A. Klain and  G.-C. Rota:}
 A continuous analogue of Sperner's theorem. 
\emph{Comm. Pure Appl. Math.} \textbf{ 50} (1997), no. 3, 205--223. 

\bibitem[Ku]{Ku} \textsc{R. K\"uhnau:} Eine Methode, die Positivit\"at einer 
Funktion zu pr\"ufen, \emph{Z. Angew. Math. Mech.} \textbf{74} (1994), no. 2, 
140--143.

\bibitem[Kum]{Kum} \textsc{E. E. Kummer:} \"Uber die hypergeometrische Reihe, 
 \emph{J. Reine Angew. Math.} \textbf{15} (1836), 39--83 and 127--172.

\bibitem[LS]{LS} \textsc{A. Laforgia and S. Sismondi:}
Some functional inequalities for complete elliptic integrals.
\emph{Rend. Circ. Mat. Palermo} (2) \textbf{41} (1992), no. 2, 302--308.

\bibitem[LV]{LV} \textsc{O. Lehto and K. I. Virtanen:} \emph{Quasiconformal 
Mappings in the Plane}, 2nd ed., Springer-Verlag, New York, 1973.

\bibitem[Le]{Le} \textsc{D. C. van Leijenhorst,} 
Algorithms for the approximation of $\pi$, 
\emph{Nieuw Archief Wisk.} \textbf{14} (1996), 255--274.


\bibitem[Lu]{Lu} \textsc{A. Lupas:} Some BBP-functions, 
Manuscript 2000. 


\bibitem[PV1]{PV1} \textsc{S. Ponnusamy and M. Vuorinen:} 
Asymptotic expansions 
and inequalities for hypergeometric functions, \emph{Mathematika} 
\textbf{44} (1997), 278--301.

\bibitem[PV2]{PV2} \textsc{S. Ponnusamy and M. Vuorinen:}
 Univalence and convexity
properties for confluent hypergeometric functions, \emph{
Complex Variables Theory Appl.} \textbf{36} (1998), 73-97.


\bibitem[PV3]{PV3} \textsc{S. Ponnusamy and M. Vuorinen:} 
Univalence and convexity
properties for Gaussian hypergeometric functions, \emph{
Rocky Mountain J. Math.} \textbf{31} (2001), 327--353.

\bibitem[PBM]{PBM} \textsc{A. P. Prudnikov, Yu. A. Brychkov, and 
O. I. Marichev:} \emph{Integrals and Series}, Vol. 3: More Special 
Functions, trans. from the Russian by G. G. Gould, Gordon and Breach, 
New York, 1988; see Math. Comp. \textbf{65} (1996), 1380--1384 for errata.

\bibitem[Q]{Q1} \textsc{S.-L. Qiu:} The proof of a conjecture on the
first elliptic integrals (in Chinese), 
\emph{J. Hangzhou Inst. of Elect. Eng.} \textbf{3} (1993), 29--36.

%\bibitem[Q2]{Q2} \textsc{S.-L. Qiu:} 
% Some properties of the gamma and psi functions, with applications,
%\emph{Manuscript, 2001. }

\bibitem[QS]{QS} \textsc{S.-L. Qiu and J.-M. Shen:} 
On two problems concerning means, 
\emph{J. Hangzhou Inst. of Elect. Eng.} \textbf{3} (1997), 1--7.

\bibitem[QVa]{QVa} \textsc{S.-L. Qiu and M. K. Vamanamurthy:}
Sharp estimates for complete elliptic integrals,
\emph{SIAM J. Math. Anal.} \textbf{27} (1996), 823--834.

\bibitem[QVV1]{QVV1} \textsc{S.-L. Qiu, M. K. Vamanamurthy, and M. Vuorinen:}
Some inequalities for the Hersch-Pfluger distortion function,
\emph{J. Inequalities Applic.} \textbf{4} (1999), 115--139.

\bibitem[QVV2]{QVV2} \textsc{S.-L. Qiu, M. K. Vamanamurthy, and M. Vuorinen:}
Some inequalities for the growth of elliptic integrals,
\emph{SIAM J. Math. Anal.} \textbf{29} (1998),  1224--1237.


\bibitem[QVu1]{QVu1} \textsc{S.-L. Qiu and M. Vuorinen:} Landen inequalities 
for hypergeometric functions, \emph{Nagoya Math. J.} \textbf{154} 
(1999), 31--56.

\bibitem[QVu2]{QVu2} \textsc{S.-L. Qiu and M. Vuorinen:} Infinite products 
and the normalized quotients of hypergeometric functions, 
\emph{SIAM J. Math. Anal.} \textbf{30} (1999),
1057--1075.

\bibitem[QVu3]{QVu3} \textsc{S.-L. Qiu and M. Vuorinen:} Duplication
inequalities for the ratios of hypergeometric functions,
\emph{Forum Math.} \textbf{12} (2000), 109--133.

\bibitem [R1] {R1} \textsc{E. D. Rainville}: \emph{ Special Functions},
MacMillan, New York, 1960.

\bibitem[R2]{R2} \textsc{E. D. Rainville:} \emph{Intermediate Differential 
Equations}, 2nd ed., Macmillan, 1964.

\bibitem[Ra1]{Ra1} \textsc{S. Ramanujan:} \emph{The Lost Notebook and Other 
Unpublished Papers}, Introduction by G. Andrews, Springer-Verlag, 
New York, 1988.

\bibitem[Ra2]{Ra2} \textsc{S. Ramanujan:} \emph{Collected papers}, 
ed. by G. S. Hardy, P. V. Seshu Aiyar, and B. M. Wilson, 
AMS Chelsea Publ. 2000, with a commentary by B. Berndt, 357--426.

\bibitem[Sa]{Sa} \textsc{E. Salamin:} Computation of $\pi$
using arithmetic-geometric mean,
\emph{Math. Comp.} \textbf{135} (1976), 565--570.


\bibitem[S1]{S1} \textsc{J. S\'andor:} On certain inequalities for means,
\emph{J. Math. Anal. Appl.} \textbf{189} (1995), 602--606.

\bibitem[S2]{S2} \textsc{J. S\'andor:} On certain inequalities for means, II,
\emph{J. Math. Anal. Appl.} \textbf{199} (1996), 629--635.

\bibitem[S3]{S3} \textsc{J. S\'andor:} On certain inequalities for means, III,
\emph{Arch. Math. (Basel)} \textbf{76} (2001), 34--40.

\bibitem[TY]{TY} \textsc{S. R. Tims and J. A. Tyrell:} Approximate evaluation 
of Euler's constant, \emph{Math. Gaz.} \textbf{55} (1971), 65--67.

\bibitem[T]{T} \textsc{Gh. Toader:} Some mean values related to the 
arithmetic-geometric mean, \emph{J. Math. Anal. Appl.} \textbf{218} 
(1998), 358--368.

\bibitem[UK]{UK} \textsc{J. S. Ume and Y.-H. Kim:} Some mean values related 
to the quasi-arithmetic mean, \emph{J. Math. Anal. Appl.} \textbf{252} 
(2000), 167--176.

\bibitem[VV]{VV} \textsc{M. K. Vamanamurthy and M. Vuorinen:} Inequalities for 
means. \emph{J. Math. Anal. Appl.} \textbf{183} (1994), 155--166.

\bibitem[V1]{V1} \textsc{M. Vuorinen:} Geometric properties of quasiconformal 
maps and special functions, I--III, 
\emph{Bull. Soc. Sci. Lett. \L ód\'z S\'er. Rech. D\'eform.} \textbf{24} 
(1997), 7--58.

\bibitem[V2]{V2} \textsc{M. Vuorinen:} Hypergeometric functions in geometric
function theory, in \emph{Proceedings of the Special Functions
and Differential Equations}, pp. 119-126, ed. by K. R. Srinivasa,
R. Jagannathan, and G. Van der Jeugy, Allied Publishers, 
New Delhi, 1998.

\bibitem[Y]{Y} \textsc{R. M. Young:} Euler's constant, \emph{Math. Gaz.} 
\textbf{75(472)} (1991), 187--190.

\end{thebibliography}
\end{document}